\documentclass[amscd,amssymb,11pt]{amsart}

\usepackage{amssymb}
\usepackage{graphicx}
\usepackage{xcolor}
\usepackage[all]{xy}

\newtheorem{thm}{Theorem}[section]
\newtheorem{lem}[thm]{Lemma}

\theoremstyle{definition}

\theoremstyle{remark}

\newtheorem{rem}[thm]{Remark}

\numberwithin{equation}{section}

\newcommand{\thmref}[1]{Theorem~\ref{#1}}

\newcommand{\secref}[1]{\S\ref{#1}}

\newcommand{\lemref}[1]{Lemma~\ref{#1}}

\newcommand{\Hom}{\operatorname{Hom}}
\newcommand{\Rep}{\operatorname{Rep}}
\newcommand{\Fix}{\operatorname{Fix}}

\newcommand{\Z}{{\mathbb  Z}}
\newcommand{\F}{{\mathbb  F}}

\newcommand{\ra}{\rightarrow}
\newcommand{\xra}{\xrightarrow}

\begin{document}

\title[Chromatic Smith Theorem]{A short proof of the chromatic Smith Fixed Point Theorem.}

\author[Kuhn]{Nicholas J.~Kuhn}
\email{njk4x@virginia.edu}

\address{Department of Mathematics \\ University of Virginia \\ Charlottesville, VA 22903}


\date{December 9, 2021.}


\begin{abstract}

We give a short and much simplified proof of the main theorem of the recent study, by T. Barthel, M. Hausmann, N. Naumann, T. Nikolaus, J. Noel, and N. Stapleton, of the Balmer spectrum for $A$--equivariant stable homotopy when $A$ is a finite abelian $p$--group.  This theorem says that if $A$ is a finite abelian $p$-group of rank $r$, and $X$ is a finite $A$--space that is acyclic in the $(n+r)$th Morava $K$--theory, then its space of fixed points, $X^A$, will be acyclic in the $n$th Morava $K$--theory.  It is a chromatic homotopy version of P.~A.~Smith's classic theorem about the mod $p$ homology of the fixed points of a finite $A$--space.
\end{abstract}

\maketitle

\section{Introduction} \label{introduction}

If $G$ is a finite group, say that $G$--space $X$ is a finite $G$--space if it is a retract of a finite $G$--CW complex in the $G$--equivariant homotopy category.  We let $X^G$ denote its subspace of fixed points.

The point of this note is to give a very short and simplified proof of the main theorem of \cite{6 author}, which can be stated as follows.

\begin{thm} \label{main thm}  Let $A$ be a finite abelian $p$--group of rank $r$, and let $X$ be a finite $A$--space. If  $\widetilde K(n+r)_*(X)=0$ then $\widetilde K(n)_*(X^A)=0$.
\end{thm}

Here $K(n)_*$ denotes the $n$th Morava $K$--theory at the prime $p$.

A number of remarks are in order.

\begin{rem}

\noindent{\bf (a)} A finite space is acyclic in mod $p$ homology if and only if it is acyclic in $K(n)_*$ for all large $n$.  Thus \thmref{main thm} (specialized to $A=C_p$ and then iterated) both implies and generalizes a classical theorem of P.~A.~Smith \cite{PA Smith 1941}: if $G$ is a finite $p$--group and $X$ is a finite $G$--space, then $\widetilde H_*(X;\Z/p) = 0$ implies that $\widetilde H_*(X^G;\Z/p) = 0$.

{\bf (b)} \ \cite[Theorem 2.1]{6 author} (or, when $A=C_p$, \cite[Proposition 7.1]{balmer sanders}) is stated in terms of the geometric fixed point functor applied to a compact object in $A$--spectra, but is easily seen to be equivalent to the theorem as stated above.  (For more detail, see \cite[\S3.1]{kuhn lloyd 2020}.)

{\bf (c)} \ By iteration, the theorem for cyclic $p$--groups implies the theorem for abelian $p$--groups, however the proof seems no easier when specialized to the cyclic case.  More generally, iterated use of the cyclic case of the theorem implies the following theorem: whenever $H$ is a subgroup of a finite $p$--group $G$ and $X$ is a finite $G$--space,
$$ \widetilde K(n+r)_*(X^H)=0 \Rightarrow \widetilde K(n)_*(X^G)=0$$
where $r$ is minimal such that there exists a sequence of subgroups $ H = K_0 \lhd K_1 \lhd \dots \lhd K_r = G$ with each $K_{i-1}$ normal in $K_i$ and  $K_i/K_{i-1}$ cyclic. (Again, this is discussed in more detail in \cite{kuhn lloyd 2020}.)

{\bf (d)} \ \cite[Theorem 2.2]{kuhn lloyd 2020} goes as follows: given a finite $p$--group $G$, $H<G$, $m$, and $n$ such that, for all finite $G$--spaces $X$,
$$ \widetilde K(m)_*(X^H) = 0 \Rightarrow \widetilde K(n)_*(X^G)=0,$$
then, if $X$ is {\em any} finite $G$--space, one has
$$  \dim_{K(m)_*} K(m)_*(X^H) \geq \dim_{K(n)_*} K(n)_*(X^G).$$
Combined with \thmref{main thm} and specialized to $G=C_{p^k}$, one deduces that
$$  \dim_{K(n+1)_*} K(n+1)_*(X) \geq \dim_{K(n)_*} K(n)_*(X^{C_{p^k}}),$$
which both implies and generalizes a classic theorem of E.E.Floyd \cite{floyd TAMS 52}: if $X$ is a finite $C_p$--space, then
$$\dim_{\Z/p} H_*(X;\Z/p) \geq \dim_{\Z/p} H_*(X^{C_p};\Z/p).$$
See \cite{kuhn lloyd 2020, kuhn lloyd 2021} for some applications.

One can also use the \cite{kuhn lloyd 2020} result in contrapositive form, and look for examples that show that the constraints in remark (c) are best possible. This is the strategy implemented in \cite{kuhn lloyd 2020}, but whether these constraints are best is still an open conjecture for many pairs $H<G$.  See \cite{kuhn lloyd 2020} for what is known about this.

\end{rem}

Our proof of \thmref{main thm} makes use of work of Stapleton in \cite{stapleton}, which generalized the character theory of \cite{hkr}. Ideas from \cite{stapleton} are certainly used in \cite{6 author}, but our proof uses this work more explicitly, while avoiding all discussion of equivariant spectra, geometric fixed point functors, generalized Tate constructions, and classifying spaces of families, among other things.

After recalling what we need of Stapleton's work in \secref{stapleton char section}, we give our proof of \thmref{main thm} in \secref{main thm proof section}.

The author would like to thank William Balderrama for conversations about this work.

\section{Stapleton's partial character map} \label{stapleton char section}

\subsection{Fixed Point Functors}  The character maps of \cite{hkr} and \cite{stapleton} both involve certain fixed point constructions, which we now describe.

If $\Gamma$ is a profinite group and $G$ is a finite group, we let $\Hom(\Gamma,G)$ denote the set of continuous homomorphisms.  $G$ acts on $\Hom(\Gamma,G)$ via conjugation, and we define a functor
$$ \Fix_{\Gamma,G}: G\text{--spaces} \ra G\text{--spaces}$$
by letting $\Fix_G(X)$ be a sub-$G$--space of $\Hom(\Gamma,G) \times X$ as follows:
$$ \Fix_{\Gamma,G}(X) = \{(\alpha,x) \ | \ \alpha(\gamma)x=x \text{ for all } \gamma \in \Gamma\} \subset \Hom(\Gamma,G) \times X.$$

As described in \cite[\S6]{kuhn char paper}, these functors satisfy various nice properties; in particular, they are homotopy functors that preserve $G$--equivariant pushouts.  Consequently, if $h^*$ is a generalized cohomology theory then
$$ X \mapsto h^*(EG \times_G \Fix_{\Gamma,G}(X))$$
will be an $G$--equivariant cohomology theory on finite $G$--spaces.

It is useful to write this a bit differently.  Let $\Rep(\Gamma,G) = \Hom(\Gamma,G)/G$. Then there is a natural homeomorphism of $G$--spaces
$$ \Fix_{\Gamma, G}(X) = \coprod_{\alpha \in \Rep(\Gamma, G)} G \times_{C_G(\alpha(\Gamma))} X^{\alpha(\Gamma)},$$
and thus a natural homeomorphism
$$ EG \times_G \Fix_{\Gamma,G}(X) = \coprod_{\alpha \in \Rep(\Gamma, G)} E C_G(\alpha(\Gamma)) \times_{C_G(\alpha(\Gamma))} X^{\alpha(\Gamma)},$$
and finally a natural decomposition
$$ h^*(EG \times_G \Fix_{\Gamma,G}(X)) = \prod_{\alpha \in \Rep(\Gamma, G)} h^*(E C_G(\alpha(\Gamma)) \times_{C_G(\alpha(\Gamma))} X^{\alpha(\Gamma)}).$$

Now suppose that $\alpha(\Gamma)=G$. Then the corresponding factor in the product is $h^*(BC_G(G) \times X^G)$.
We specialize this observation to the case when $\Gamma = \Z_p^r$, the product of $r$ copies of the $p$-adic integers $\Z_p$. Since there exists a continuous surjective homomorphism $\alpha: \Z_p^r \twoheadrightarrow A$ if and only if $A$ is an abelian $p$--group of rank at most $r$,  we deduce the following.

\begin{lem} \label{factor lemma} If $A$ is an abelian $p$--group of rank at most $r$, and $h^*$ is a generalized cohomology theory, then $h^*(BA \times X^A)$ is naturally a direct summand in $h^*(EA \times_A \Fix_{\Z_p^r,A}(X))$ for all $A$--spaces $X$.
\end{lem}

\subsection{Transchromatic characters} \label{stapleton subsection}

We describe what we need from \cite{stapleton}.   It is convenient to fix $n$ and $r$, and then let $E=E_{n+r}$, the Morava $E$--theory of height $(n+r)$, and $LE=L_{K(n)}E_{n+r}$, the Bousfield localization of $E_{n+r}$ with respect to the $r$th Morava $K$--theory.

It is illuminating to recall (\cite[\S2.1]{stapleton}) that
$$ E^0 = \mathbb W(\F_{p^{n+r}})[[u_1, \dots, u_{n+r}]]$$
and (\cite[Theorem 1.5.4]{hovey}) that
$$ LE^0 = \mathbb W(\F_{p^{n+r}})[[u_1, \dots, u_{n+r}]][u_n^{-1}]^{\wedge}_{I_n},$$
where $I_n = (p,u_0, \dots, u_{n-1})$.

\begin{lem} \label{LE(BA) lem} $LE$ is complex oriented of height $n$. If $A$ is a finite abelian $p$--group, then $LE^*(BA)$ is a finitely generated free $LE^*$--module, and thus there is a Kunneth ismorphism
$$ LE^*(BA) \otimes_{LE^*}LE^*(Y) \simeq LE^*(BA \times Y)$$
for all spaces $Y$.
\end{lem}

The first of these properties is clear, and the second is well known (\cite[Proposition 5.2]{hkr}) and implies the last statement.

Stapleton constructs a faithfully flat $LE^0$--algebra $C$ and proves the following theorem.

\begin{thm} \label{stapleton thm} \cite[Main Theorem]{stapleton}  For any finite group $G$ and finite $G$--space $X$, there is a natural algebra isomorphism
$$ \chi: E^*(EG \times_G X)\otimes_{E^0}C \xra{\sim} LE^*(EG \times_G \Fix_{\Z_p^r, G}(X))\otimes_{LE^0}C.$$
\end{thm}

The case $n=0$ (when $LE$ is a rational theory) was shown in \cite{hkr}.

\section{Proof of \thmref{main thm}} \label{main thm proof section}

Let $A$ be a finite abelian $p$--group of rank $r$, and let $X$ be a finite $A$--space. Our goal is to show that if  $\widetilde K(n+r)_*(X)=0$ then $\widetilde K(n)_*(X^A)=0$.

As in \secref{stapleton subsection}, we let $E=E_{n+r}$, $LE = L_{K(n)}E_{n+r}$, and $C$ be Stapleton's flat $LE^0$--algebra as in \thmref{stapleton thm}.

Our assumption that $\widetilde K(n+r)_*(X)=0$ is equivalent to saying that the constant map $X \ra *$ induces an isomorphism
$$ K(n+r)_*(X) \xra{\sim} K(n+r)_*.$$

This then implies that the induced map $EA \times_A X \ra BA$ induces an isomorphism
$$ K(n+r)_*(EA \times_A X) \xra{\sim} K(n+r)_*(BA).$$
This is an immediate consequence of the fact that, given any finite group $G$, $G$--space $X$, and generalized homology theory $h_*$, there is a natural spectral sequence converging to $h_*(EG \times_G X)$ with $E^2_{*,*} = H_*(BG;h_*(X))$.

As $K(n+r)_*$--isomorphisms are $E^*$--isomorphisms, we conclude that the map $EA \times_A X \ra BA$ induces an isomorphism
$$ E^*(BA) \xra{\sim} E^*(EA \times_A X).$$

Now consider the commutative diagram
\begin{equation*}
\SelectTips{cm}{}
\xymatrix{
E^*(BA) \otimes_{E^0} C \ar[d]^{\chi} \ar[r] & E^*(EA \times_A X) \ar[d]^{\chi} \otimes_{E^0} C  \\
LE^*(EA \times_A \Fix_{\Z_p^r,A}(*)) \otimes_{LE^0}C \ar[r] & LE^*(EA \times_A \Fix_{\Z_p^r,A}(X))\otimes_{LE^0}C. }
\end{equation*}
We have just seen that the top map is an isomorphism.  So are the two vertical maps by \thmref{stapleton thm}.  Thus the bottom map is also an isomorphism.

Thanks to \lemref{factor lemma}, we can conclude that there is an isomorphism
$$ LE^*(BA) \otimes_{LE^0}C \xra{\sim} LE^*(BA \times X^A) \otimes_{LE^0}C$$
and, thanks to \lemref{LE(BA) lem}, this can be rewritten as an isomorphism
$$ LE^*(BA) \otimes_{LE^0}C \xra{\sim} LE^*(BA) \otimes_{LE^0} LE^*(X^A) \otimes_{LE^0}C.$$

Applying $LE^* \otimes_{LE^*(BA)} \text{\underline{\hspace{.1in}}}$ to this, we deduce that there is an isomorphism
$$ LE^* \otimes_{LE^0}C \xra{\sim} LE^*(X^A) \otimes_{LE^0}C, $$
and since $C$ is a faithful $LE^0$--module, we further deduce that there is an isomorphism
$$ LE^* \xra{\sim} LE^*(X^A). $$
Finally, since $LE^*$ has height $n$, an $LE^*$--isomorphism will be a $K(n)_*$--isomorphism, so we conclude that $K(n)_*(X^A) \ra K(n)_*$ is an isomorphism, and thus $\widetilde K(n)_*(X^A) = 0$.


\begin{thebibliography}{99}



\bibitem[BS17]{balmer sanders} P.Balmer and B.Sanders, {\em The spectrum of the equivariant stable homotopy theory of a finite group}, Invent. Math. {\bf 208} (2017), 283--326.

\bibitem[6A19]{6 author} T. Barthel, M. Hausmann, N. Naumann, T. Nikolaus, J. Noel, and N. Stapleton, {\em The Balmer spectrum of the equivariant homotopy category of a finite abelian group}, Invent. Math. {\bf 216} (2019), 215-240.

\bibitem[F52]{floyd TAMS 52} E.~E.~Floyd, {\em On periodic maps and the Euler characteristics of associated spaces}, Trans. Amer. Math. Soc. {\bf 72} (1952), 138-147.

\bibitem[HKR00]{hkr} M.~J.~ Hopkins, N.~J.~ Kuhn, and D.~C.~ Ravenel, {\em Generalized group characters and complex oriented cohomology theories},  J. Amer. Math. Soc. {\bf 13} (2000), 553-594.

\bibitem[H97]{hovey}  M.~A.~Hovey, {\em $v_n$-elements in ring spectra and applications to bordism theory}, Duke Math. J. {\bf 88} (1997), 327-356.

\bibitem[K89]{kuhn char paper} N.~J.~Kuhn, {\em Character rings in algebraic topology}, in Advances in Homotopy Theory (Cortona, 1988), 111-126,
London Math. Soc. Lecture Note Ser., 139, Cambridge Univ. Press, Cambridge, 1989.

\bibitem[KL20]{kuhn lloyd 2020} N.~J.~Kuhn and C.~J.~R.~Lloyd, {\em Chromatic Fixed Point Theory and the Balmer spectrum for extraspecial 2-groups}, preprint 2020.

\bibitem[KL21]{kuhn lloyd 2021} N.~J.~Kuhn and C.~J.~R.~Lloyd, {\em Computing the Morava K-theory of real Grassmanians using chromatic fixed point theory}, preprint 2021.

\bibitem[S41]{PA Smith 1941} P. ~A. ~Smith, {\em Fixed point theorems for periodic transformations}, Amer.J.Math. {\bf 63}(1941), 1--8.

\bibitem[S13]{stapleton}  N.~ Stapleton, {\em Transchromatic generalized character maps},  Algebr. Geom. Topol. {\bf 13} (2013), 171-203.


\end{thebibliography}
\end{document}